\documentclass[11pt]{article}
\usepackage{amsmath,latexsym,epsf,mathabx}
\usepackage[all]{xy}
\usepackage{mathrsfs}
\usepackage{soul}
\usepackage{color}
\begin{document}

\newcommand{\nc}{\newcommand}
\newtheorem{Th}{Theorem}[section]
\newtheorem{Def}[Th]{Definition}
\newtheorem{Lem}[Th]{Lemma}
\newtheorem{Pro}[Th]{Proposition}
\newtheorem{Cor}[Th]{Corollary}
\newtheorem{Rem}[Th]{Remark}
\newtheorem{Exm}[Th]{Example}
\newtheorem{Sc}[Th]{}
\def\Pf#1{{\noindent\bf Proof}.\setcounter{equation}{0}}
\def\bskip#1{{ \vskip 20pt }\setcounter{equation}{0}}
\def\sskip#1{{ \vskip 5pt }\setcounter{equation}{0}}
\def\mskip#1{{ \vskip 10pt }\setcounter{equation}{0}}
\def\bg#1{\begin{#1}\setcounter{equation}{0}}
\def\ed#1{\end{#1}\setcounter{equation}{0}}

\soulregister\cite7 
\soulregister\citep7 
\soulregister\citet7 
\soulregister\ref7


\title{\bf  Gorenstein projective objects and recollements of Abelian categories
\thanks{Supported by the Top talent project of AHPU in 2020, the National Natural Science Foundation of China (Grants No.11801004)
and the Startup Foundation for Introducing Talent of AHPU (Grant: 2020YQQ067)
. }
}
\smallskip
\author{ Peiyu Zhang\thanks{Correspondent author}, Qianqian Shu and Dajun Liu\\ 
\footnotesize ~E-mail:~zhangpy@ahpu.edu.cn, 15955660648@163.com and liudajun@ahpu.edu.cn; \\
\footnotesize School of Mathematics and Physics,  Anhui Polytechnic University, Wuhu, China.
}

\date{}
\maketitle
\baselineskip 15pt
%
%
\begin{abstract}
\vskip 5pt%
In this paper, we study the relationship of Gorenstein projective objects among three
Abelian categories in a recollement. As an application, we introduce the relation of $n$-Gorenstein tilting modules
(and Gorenstein syzygy modules) in three Abelian categories. For a recollement $R(\mathscr A', ~\mathscr A, ~\mathscr A'')$ of Abelian categories,
we show that a resolving subcategory in $\mathscr A$ induce
two resolving subcategories in $\mathscr A'$ and $\mathscr A''$. On the other hand, we also prove that two resolving subcategories in $\mathscr A'$ and
$\mathscr A''$ can induce a resolving subcategory in $\mathscr A$.
Moreover, we give the size relationship between the relative global dimensions of three Abelian categories.

\mskip\




\end{abstract}
%
\vskip 30pt

\section{Introduction}

Gorenstein homological algebra is the relative version of homological
algebra that uses Gorenstein projective, Gorenstein injective and Gorenstein
flat modules instead of the classical projective, injective and flat modules.
Enochs and Jenda introduced in \cite{EJ1} Gorenstein projective modules for arbitrary modules over a general
ring, which is a generalization of finitely generated modules of Gorenstein dimension zero.
The Gorenstein projective objects were studied by many authors.
For example, Holm gave some homological properties of Gorenstein projective, and proved that
every $R$-module $M$ with finite Gorenstein projective dimensions admits a surjective Gorenstein projective modules precover
in \cite{Holm}. Using the relative homological theory developed by Enochs and Jenda \cite{EJ}, the authors \cite{YLO} introduced Gorenstein
cotilting and tilting modules, and gave a characterization of Gorenstein tilting modules which similar to Bazzoni characterization
of $n$-tilting modules \cite{BS}, and so on.

A recollement of Abelian categories is a kind of algebraic structure consisting of three Abelian categories and six functors.
It has applications in many aspects. Recollements of Abelian categories were studied by many authors, see \cite{FP,LYN,mh,PS,PC,PC1}.
In this paper, all Abelian categories has enough projective and injective objects. In section 3, we mainly consider the relationship of Gorenstein projective objects
($n$-Gorenstein tilting modules) in a recollement and give the following conclusion. In section 4, we give two applications on Theorem \ref{Th11}.

\bg{Th}\label{Th11}%
In a recollement $R(\mathscr A', ~\mathscr A, ~\mathscr A'')$, let $X'$, $X$ and $X''$ be a Gorenstein projective objects
in $\mathscr A'$, $\mathscr A$ and $\mathscr A''$, respectively.

$(1)$ If $j_{\ast}$ is exact, then $j^{\ast}X$ is a Gorenstein projective object in $\mathscr A''$;

$(2)$ If $i^{\ast}$ is weak perfect and $i_{\ast}$ satisfies the condition $\mathrm{(SP2)}$, then $i^{\ast}X$ is a Gorenstein projective object in $\mathscr A'$.

$(3)$ If $i^{!}$ is exact and $i^{\ast}$ is fully faithful, then $i_{\ast}X'$ is a Gorenstein projective object in $\mathscr A$.

$(4)$ If $j_{!}$ is weak perfect and $j^{\ast}$ satisfies the condition $\mathrm{(SP2)}$, then $j_{!}X''$ is a Gorenstein projective object in $\mathscr A$.
\ed{Th}

In \cite{Holm}, the author proved that Gorenstein projective objects has the following properties: (1) All projective objects are Gorenstein projective;
(2) Gorenstein projective objects is closed under direct summands; (3) for any exact sequence $\xymatrix{0\ar[r]&A\ar[r] &B\ar[r] &C\ar[r] &0}$ with $C$ being Gorenstein projective,
then $A$ is Gorenstein projective if and only if $B$ is Gorenstein projective.
In \cite{AR}, an subcategory is said to be a resolving subcategory, if it satisfies the above three conditions.
Clearly, the subcategory consist of all Gorenstein projective objects is a resolving subcategory. In section 4, we mainly prove the following two results.

\bg{Th}\label{Th12}%
Let $R(\mathscr A', ~\mathscr A, ~\mathscr A'')$ be a recollement of Abelian categories.

$(1)$ If $\mathscr X$ be a resolving subcategory of $\mathscr A$ and $i_{\ast}i^{\ast}\mathscr X\subseteq \mathscr X$, then $i^{\ast}\mathscr X$ is resolving in $\mathscr A'$.

$(2)$ If $\mathscr X$ be a resolving subcategory of $\mathscr A$, $j_{\ast}$ is exact and $j_{!}j^{\ast}\mathscr X\subseteq \mathscr X$ or $j_{\ast}j^{\ast}\mathscr X\subseteq \mathscr X$,
then $j^{\ast}\mathscr X$ is resolving in $\mathscr A''$.

$(3)$ Assume that $\mathscr X'$ and $\mathscr X''$ are resolving in $\mathscr A'$ and $\mathscr A''$, respectively. If both $i^{\ast}$ and $j_{\ast}$ are exact, then the following subcategory

\centerline{$\mathscr X:=\{X\in \mathscr A|~i^{\ast}X\in \mathscr X'$ and $j^{\ast}X\in \mathscr X''\}$}
\noindent is a resolving subcategory of $\mathscr A$.
\ed{Th}

\bg{Th}\label{th13}%
Let $R(\mathscr A', ~\mathscr A, ~\mathscr A'')$ be a recollement of Abelian categories and
$\mathscr X$ be a resolving subcategory of $\mathscr A$.

$(1)$ If $i_{\ast}i^{\ast}\mathscr X\subseteq \mathscr X$,
$j_{!}j^{\ast}\mathscr X\subseteq \mathscr X$ and $j_{\ast}$ is exact, then the following inequations hold.

$~~~~$$(i)$ $\mathrm{gl. dim}_{\mathscr X}\mathscr A\leq \mathrm{gl. dim}_{i^{\ast}\mathscr X}\mathscr A'+\mathrm{gl. dim}_{j^{\ast}\mathscr X}\mathscr A''+1.$

$~~~~$$(ii)$ $\mathrm{gl. dim}_{j^{\ast}\mathscr X}\mathscr A''\leq \mathrm{gl. dim}_{\mathscr X}\mathscr A$.

$(2)$ If $i_{\ast}i^{\ast}\mathscr X\subseteq \mathscr X$ and $i^{\ast}$ is exact,
then $\mathrm{gl. dim}_{i^{\ast}\mathscr X}\mathscr A'\leq\mathrm{gl. dim}_{\mathscr X}\mathscr A$.

\ed{Th}

\section{Gorenstein projective objects}

Given two functors $F$: ${\mathscr A}\longrightarrow {\mathscr B}$ and the functor $G$: ${\mathscr B}\longrightarrow {\mathscr A}$, where both $\mathscr A$ and $\mathscr B$
are Abelian categories. we said $(F$, $G)$ to be an adjoint pair,
if there is a natural isomorphism $\sigma_{X,Y}$: $\mathrm{Hom}_{{\mathscr B}}(FX, Y)\cong\mathrm{Hom}_{{\mathscr A}}(X, GY)$ for any $X\in {\mathscr A}$ and $Y\in {\mathscr B}$.
In this case, if there is a functor $H$: ${\mathscr A}\longrightarrow {\mathscr B}$ such that the pair $(G,~H)$ is also an adjoint pair, then we call $(F,~G,~H)$ being adjoint
triple as follows.
$$\xymatrix{
{\mathscr B} \ar[rrr]|{\ G}&&& {\mathscr A}
\ar@/^/@<2ex>[lll]|{H}\ar@/_/@<-2ex>[lll]|{F} }$$

\mskip\
For an Abelian category $\mathscr A$, we denotes all projective objects by $\mathrm{Proj}(\mathscr A)$.
Let $\mathscr C$ is an subcategory of $\mathscr A$, the class $\mathrm{add} \mathscr C$ is all direct summands of $C^{n}$, where $C\in \mathscr C$ and $n$ is a positive integer.
Firstly, we give the following statements on adjoint triple, which is very useful in the rest of the section.

\bg{Lem}\label{lem1}
Assume that $(F,~G,~H)$ is an adjoint triple, then the following conclusions hold.

$(1)$ If $H$ is an exact functor, then $G$ preserves projective objects.
Moreover, if $F$ is fully faithful and $\mathscr B$ has enough projective object,
then $\mathrm{Proj}(\mathscr A)=\mathrm{add}[G(\mathrm{Proj}(\mathscr B))]$.

$(2)$ The functor $F$ preserves projective objects. If $G$ is fully faithful and $\mathscr A$ has enough projective objects,
then $\mathrm{Proj}(\mathscr B)=\mathrm{add}[F(\mathrm{Proj}(\mathscr A))]$.

\ed{Lem}

\Pf. (1) For any exact sequence $0\longrightarrow X\longrightarrow Y \longrightarrow Z\longrightarrow0$ in $\mathscr A$ and any projective object $P\in\mathscr B$.
We consider the following communicative diagram
$$\xymatrix{
0\ar[r]&\mathrm{Hom}_{\mathscr A}(GP,~X)\ar[r]\ar^{\cong}[d]&\mathrm{Hom}_{\mathscr A}(GP,~Y)\ar[r]\ar^{\cong}[d]&\mathrm{Hom}_{\mathscr A}(GP,~Z)\ar@{.>}[r]\ar^{\cong}[d]&0\\
0\ar[r]&\mathrm{Hom}_{\mathscr B}(P,~HX)\ar[r]&\mathrm{Hom}_{\mathscr B}(P,~HY)\ar[r]&\mathrm{Hom}_{\mathscr B}(P,~HZ)\ar[r]&0
}$$
Since $H$ is an exact functor, the second row is exact. Thus, the first row is also exact. i.e.,
$GP$ is a projective object in $\mathscr A$.

For any projective object $P$ in $\mathscr A$, there is a surjective morphism $Q\longrightarrow FP\longrightarrow0$ with $Q\in\mathrm{Proj}(\mathscr B)$.
And then we can obtain an exact sequence $GQ\longrightarrow GFP\longrightarrow0$ since $G$ is exact. Since $F$ is fully faithful, $GFP\cong P$ is projective,
thus $P$ is a direct summand of $GQ$. i.e., $\mathrm{Proj}(\mathscr A)=\mathrm{add}[G(\mathrm{Proj}(\mathscr B))]$.

(2) The proof similar to (1).
\ \hfill $\Box$

%
%
%

\mskip\

Recall that a complete projective resolution is an exact sequence of projective objects,
$P^{\bullet}=:$ $ \cdots \longrightarrow P_{1}\longrightarrow P_{0}\longrightarrow P_{-1}\longrightarrow \cdots\ $,
such that $\mathrm{Hom}(P^{\bullet},~Q)$ is exact for any projective module $Q$.
An object $M$ is called Gorenstein projective \cite{Holm}, if there is a complete projective resolution $P^{\bullet}$
with $M\cong \mathrm{Im}(P_{0}\longrightarrow P_{-1})$. Dually, we can define complete injective resolution
and Gorenstein injective object.

Motivated by the definition of compatible bimodules in \cite[Definition 1.1]{ZP}, Peng, Zhu and Huang introduced the perfect functor in \cite[Definition 3.3]{PZH}.
i.e., The functor $F: \mathscr A\longrightarrow \mathscr B$ is called perfect if the following two conditions are satisfied.
(P1) If $Q^{\bullet}$ is an exact sequence of projective objects, then $FQ^{\bullet}$ is exact;
(P2) If $P^{\bullet}$ is a complete $\mathscr B$-projective resolution, then $\mathrm{Hom}_{\mathscr B}(P^{\bullet},FQ)$ is exact for any $Q$ being projective object of $\mathscr A$.
If the functor $F$ satisfies the condition (P1) , then we said to be weak perfect. Here we give the definition of strong perfect functors as follows.

\bg{Def}\label{storngf}
The functor $F: \mathscr A\longrightarrow \mathscr B$ is called strong perfect if the following two conditions are satisfied.

$\mathrm{(SP1)}$ $F$ is an exact functor.

$\mathrm{(SP2)}$ F satisfies the condition $\mathrm{(P2)}$.
\ed{Def}

\bg{Rem}\label{storngf}
$(1)$ By the related definitions, it is easy to see that the strong perfect functor is perfect.

$(2)$ Let $\Lambda$ and $\Gamma$ be artin algebras, and $M$ be a compatible $(\Lambda,~\Gamma)$-bimodule, then $M\otimes_{\Gamma}-$ is perfect \cite{PZH}.
Furthermore, if $M$ is flat $\Gamma$-module, then $M\otimes_{\Gamma}$- is strong perfect.
\ed{Rem}

\bg{Lem}\label{lem2}
Assume that $(F,~G,~H)$ is an adjoint triple. If $H$ is an exact functor and $F$ is fully faithful, then
$H$ satisfies the condition $(\mathrm{SP2})$, i.e., $H$ is strong perfect.
\ed{Lem}

\Pf. For any complete projective resolution $\xymatrix{\cdots\ar[r]&P_{-1}\ar[r] &P_{0}\ar[r]^{f_{0}}&P_{1}\ar[r] &\cdots}$ in $\mathscr B$,
where $X_{i}=\mathrm{Ker}f_{i}$ for any $i$. For any $Q\in \mathrm{Proj}(A)$ such that there is a projective object $P$ satisfying $GP\cong Q$,
we consider the following commutative diagram with exact row.
$$\xymatrix{
0\ar[r]&\mathrm{Hom}_{\mathscr B}(X_{1},~HQ)\ar[d]^{\cong}\ar[r]&\mathrm{Hom}_{\mathscr B}(P_{0},~HQ)\ar[d]^{\cong}\ar[r]&\mathrm{Hom}_{\mathscr B}(X_{0},~HQ)\ar[d]^{\cong}\\
0\ar[r]&\mathrm{Hom}_{\mathscr A}(GX_{1},~Q)\ar[r]\ar[d]^{\cong}&\mathrm{Hom}_{\mathscr A}(GP_{0},~Q)\ar[r]\ar[d]^{\cong}&\mathrm{Hom}_{\mathscr A}(GX_{0},~Q)\ar[d]^{\cong}\\
0\ar[r]&\mathrm{Hom}_{\mathscr A}(GX_{1},~GP)\ar[r]&\mathrm{Hom}_{\mathscr A}(GP_{0},~GP)\ar[r]&\mathrm{Hom}_{\mathscr A}(GX_{0},~GP)
}$$
Noth that $\xymatrix{0\ar[r]&\mathrm{Hom}_{\mathscr B}(X_{1},~P)\ar[r]&\mathrm{Hom}_{\mathscr B}(P_{0},~P)\ar[r]&\mathrm{Hom}_{\mathscr B}(X_{0},~P)\ar[r]&0}$ is a short
exact sequence and $G$ is exact. So the third row in above diagram is a short exact sequence, and then the first row in above diagram is also a short exact sequence.
Consequently, $H$ satisfies the condition $(\mathrm{SP2})$.

\bg{Pro}\label{pro1}
Assume that $(F,~G,~H)$ is an adjoint triple and $X$ (resp., $Y$) is a Gorenstein projective object in $\mathscr B$ (resp., $\mathscr A$).
Then we have that the following conclusions hold.

$(1)$ If $H$ is an exact functor and $F$ is fully faithful, then GX is a Gorenstein projective object.

$(2)$ If F is weak perfect and G satisfies the condition $(\mathrm{SP2})$, then $FY$ is a Gorenstein projective object in $\mathscr B$.
\ed{Pro}

\Pf. (1) Since $X$ is Gorenstein projective, there is a complete $\mathscr B$-projective resolution
$$\xymatrix{
\cdots\ar[r]&P_{-1}\ar[r]^{f_{-1}} &P_{0}\ar[r]^{f_{0}} &P_{1}\ar[r]^{f_{1}} &\cdots\ ,
}$$
where $\mathrm{Ker}f_{i}=X_{i}$ for any $i$ and $X_{0}=X$. Consider the exact sequence
$0\longrightarrow X_{0}\longrightarrow P_{0} \longrightarrow X_{1}\longrightarrow0$.
Note $G$ is exact and preserves projective objects by Lemma \ref{lem1}. We can obtain a new exact sequence $0\longrightarrow GX_{0}\longrightarrow GP_{0} \longrightarrow GX_{1}\longrightarrow0$
with $GP_{0}$ being projective in $\mathscr A$. For any projective object $P\in\mathscr A$, we have the following commutative diagram:
$$\xymatrix{
0\ar[r]&\mathrm{Hom}_{\mathscr B}(GX_{1},~P)\ar[d]^{\cong}\ar[r]&\mathrm{Hom}_{\mathscr B}(GP_{0},~P)\ar[d]^{\cong}\ar[r]&\mathrm{Hom}_{\mathscr B}(GX_{0},~P)\ar[d]^{\cong}\\
0\ar[r]&\mathrm{Hom}_{\mathscr A}(X_{1},~HP)\ar[r]&\mathrm{Hom}_{\mathscr A}(P_{0},~HP)\ar[r]&\mathrm{Hom}_{\mathscr A}(X_{0},~HP)
}$$
By the lemma \ref{lem2}, $H$ satisfies the condition $(\mathrm{SP2})$. We know that the second row in above diagram is exact. i.e., The first row is exact.
So the following sequence induced by the complete projective resolution of $X$ is a complete projective resolution of $GX$
$$\xymatrix{
\cdots\ar[r]&GP_{-1}\ar[r] &GP_{0}\ar[r] &GP_{1}\ar[r] &\cdots\ ,
}$$
i.e., $GX$ is a Gorenstein projective object.

(2) Since $Y$ is Gorenstein projective in $\mathscr A$, there is a complete projective resolution
$$\xymatrix{
\cdots\ar[r]&Q_{-1}\ar[r]^{g_{-1}} &Q_{0}\ar[r]^{g_{0}} &Q_{1}\ar[r]^{g_{1}} &\cdots\ ,
}$$
where $\mathrm{Ker}g_{i}=Y_{i}$ for any $i$ and $Y_{0}=Y$. Consider the exact sequence
$0\longrightarrow Y_{0}\longrightarrow Q_{0} \longrightarrow Y_{1}\longrightarrow0$.
Since $F$ is weak perfect, we can obtain an exact sequence $0\longrightarrow FY_{0}\longrightarrow FQ_{0} \longrightarrow FY_{1}\longrightarrow0$
with $FQ_{0}$ projective in $\mathscr B$ by Lemma \ref{lem1}. For any projective object $Q\in\mathscr B$, we have the following commutative diagram:
$$\xymatrix{
0\ar[r]&\mathrm{Hom}_{\mathscr B}(FY_{1},~Q)\ar[d]^{\cong}\ar[r]&\mathrm{Hom}_{\mathscr B}(FQ_{0},~Q)\ar[d]^{\cong}\ar[r]&\mathrm{Hom}_{\mathscr B}(FY_{0},~Q)\ar[d]^{\cong}\ar@{.>}[r]&0\\
0\ar[r]&\mathrm{Hom}_{\mathscr A}(Y_{1},~GQ)\ar[r]&\mathrm{Hom}_{\mathscr A}(Q_{0},~GQ)\ar[r]&\mathrm{Hom}_{\mathscr A}(Y_{0},~GQ)\ar[r]&0
}$$
Since $G$ satisfies the condition $(\mathrm{SP2})$, the second row is exact. i.e., The first row is exact.
So the following sequence is a complete projective resolution of $FY$
$$\xymatrix{
\cdots\ar[r]&FQ_{-1}\ar[r] &FQ_{0}\ar[r] &FQ_{1}\ar[r] &\cdots\ ,
}$$
i.e., $FY$ is a Gorenstein projective.
\ \hfill $\Box$

\mskip\

In the (2) of proposition above, if we replace the condition "G satisfies the condition (SP2)" by the condition
"$H$ is exact", then we can also prove that $FY$ is Gorenstein projective similarly.

\bg{Def}\rm(\cite{FP,PC1})\label{RE}%
A recollement of an abelian category ${\mathscr A}$ by abelian categories ${\mathscr A'}$ and
${\mathscr A''}$, denoted by $R(\mathscr A', ~\mathscr A, ~\mathscr A'')$, is a diagram of additive functors as follows, satisfies the conditions below.

$$\xymatrix{
{\mathscr A^\prime} \ar[rrr]|{\ i_{\ast}}&&& {\mathscr A}
\ar@/^/@<2ex>[lll]|{i^{!}}\ar[rrr]|{\
j^{\ast}}\ar@/_/@<-2ex>[lll]|{i^{\ast}} &&& {\mathscr A''}
\ar@/_/@<-2ex>[lll]|{j_{!}}\ar@/^/@<2ex>[lll]|{j_{\ast}} }$$

$(i)$ ($i^{\ast}$, $i_{\ast}$, $i^{!}$) and ($j_{!}$, $j^{\ast}$, $j_{\ast}$) are adjoint triples;

$(ii)$ The functors $i_{\ast}$, $j_{!}$, and $j_{\ast}$ are fully faithful;

$(iii)$ $\mathrm{Im}i_{\ast}$ = $\mathrm{Ker}j^{\ast}$.

\ed{Def}

\mskip\
There are many examples of recollements, see \cite{PC}.
Next, we collect some properties of recollements, which is very useful, see \cite{FP,mh,PC,PSS,PC1}.

\bg{Pro}\label{propR}%

Let $R(\mathscr A', ~\mathscr A, ~\mathscr A'')$ be a recollement of abelian categories. Then we have the following properties.

$(1)$ $i^{\ast}j_{!}=0$ and $i^{!}j_{\ast}=0$;

$(2)$  $i^{\ast}$ and $j_{!}$ are right exact, $i^{!}$ and $j_{\ast}$ are left exact, $i_{\ast}$ and $j^{\ast}$ are exact;

$(3)$ These natural transformations $i^{\ast}i_{\ast}\longrightarrow \mathrm{Id}_{\mathscr A'}$, $\mathrm{Id}_{\mathscr A'}\longrightarrow i^{!}i_{\ast}$,
$j^{\ast}j_{\ast}\longrightarrow \mathrm{Id}_{\mathscr A''}$, $\mathrm{Id}_{\mathscr A''}\longrightarrow j^{\ast}j_{!}$ are natural isomorphisms.

$(4)$ If $i^{\ast}$ is exact, then $i^{!}j_{!}=0$ and $j_{!}$ is exact; If $i^{!}$ is exact, then $i^{\ast}j_{\ast}=0$ and $j_{\ast}$ is exact.

$(5)$ for any $A\in\mathscr A$, there exist two exact sequences
$$\xymatrix{0\ar[r]&i_{\ast}A'\ar[r]&j_{!}j^{\ast}A\ar[r]&A\ar[r]&i_{\ast}i^{\ast}A\ar[r]&0}$$
$$\xymatrix{0\ar[r]&i_{\ast}i^{!}A\ar[r]&A\ar[r]&j_{\ast}j^{\ast}A\ar[r]&i_{\ast}A_{1}'\ar[r]&0}$$
where $A_{1}'$ and $A'$ are in $\mathscr A'$.
\ed{Pro}

By Proposition \ref{pro1} and \ref{propR}, we obtain the following main result immediately.

\bg{Th}\label{result1}%
In a recollement $R(\mathscr A', ~\mathscr A, ~\mathscr A'')$, let $X'$, $X$ and $X''$ be a Gorenstein projective objects
in $\mathscr A'$, $\mathscr A$ and $\mathscr A''$, respectively.

$(1)$ If $j_{\ast}$ is exact, then $j^{\ast}X$ is a Gorenstein projective object in $\mathscr A''$;

$(2)$ If $i^{\ast}$ is weak perfect and $i_{\ast}$ satisfies the condition $\mathrm{(SP2)}$, then $i^{\ast}X$ is a Gorenstein projective object in $\mathscr A'$.

$(3)$ If $i^{!}$ is exact and $i^{\ast}$ is fully faithful, then $i_{\ast}X'$ is a Gorenstein projective object in $\mathscr A$.

$(4)$ If $j_{!}$ is weak perfect and $j^{\ast}$ satisfies the condition $\mathrm{(SP2)}$, then $j_{!}X''$ is a Gorenstein projective object in $\mathscr A$.
\ed{Th}

%
%

\section{Application}

In this section, $\mathscr A', ~\mathscr A$ and $\mathscr A''$ are module categories.

\noindent\textbf{3.1 Gorenstein tilting modules}
Replace projective modules by Gorenstein projective modules, we can compute right derived functors $\mathrm{Gext}^{i}(-,~-)$ for any $i>0$, see Chapter 12 in \cite{EJ}.
Denoted all Gorenstein projective modules (resp. Gorenstein injective modules) by $\mathcal{GP}$ (resp. $\mathcal{GI}$).
If $\mathrm{Gext}^{i}(T,~M)=0$ for any $i>0$, then we denote $M\in T^{G\bot}$.
Recall that a short exact sequence $0\longrightarrow X\longrightarrow Y\longrightarrow Z\longrightarrow0$ is called $G$-exact if it is in $\mathrm{Gext}^{1}(Z,X)$.

\mskip\

The following lemma is a characterization of $G$-exact sequences.

\bg{Lem}\label{G-exact}
The following statements are equivalent for an exact sequence $0\longrightarrow X\stackrel{i}{\longrightarrow} Y\stackrel{\pi}{\longrightarrow} Z\longrightarrow 0$:

$1)$ The sequence is G-exact;

$2)$ $\xymatrix{0\ar[r]&(P,X)\ar[r]&(P,Y)\ar[r]&(P,Z)\ar[r]&0}$ is exact for all $P\in \mathcal{GP}$;

$3)$ $\xymatrix{0\ar[r]&(Z,I)\ar[r]&(Y,I)\ar[r]&(X,I)\ar[r]&0}$ is exact for all $I\in \mathcal{GI}$;
\ed{Lem}

\Pf. By Proposition 1.5 in \cite{AS} and Lemma 2.1 in \cite{zpy0}.
\ \hfill $\Box$

\bg{Def}\rm(\cite[Definition ~3.2]{YLO})\label{gtilting}
An R-module T is called an $n$-Gorenstein tilting module if it
satisfies the following three conditions:

$i)$ $\mathrm{pd}_{G}T \leq n$;

$ii)$ $\mathrm{Gext}^{i}(T,T^{(I)})=0$ for each $i>0$ and all sets I;

$iii)$ There exists a long G-exact sequence $0\longrightarrow P\longrightarrow T^{0}\longrightarrow T^{1}\longrightarrow\cdots \longrightarrow T^{n}\longrightarrow 0$ with
each $T^{i}\in\mathrm{Add}T$ and any $P\in \mathcal{GP}$ .

If T satisfies the conditions $i)$ and $ii)$, then we say T to be partial $n$-Gorenstein tilting.
\ed{Def}

\bg{Th}\label{result2}%
In a recollement $R(\mathscr A', ~\mathscr A, ~\mathscr A'')$,
let $T\in\mathscr A$ be an $n$-Gorenstein tilting module.

$(1)$ If $i^{\ast}$ is exact, $i_{\ast}$ satisfies the condition $\mathrm{(SP2)}$ and $i_{\ast}i^{\ast}T\subseteq T^{G\bot}$, then $i^{\ast}T$ is partial $n$-Gorenstein tilting;

$(2)$ If $j_{\ast}$ is exact, $j_{!}$ is weak perfect and $j_{\ast}j^{\ast}T\subseteq T^{G\bot}$, then $j^{\ast}T$ is $n$-Gorenstein tilting;
\ed{Th}

\Pf. (1) Since $\mathrm{pd}_{G}T \leq n$, there is a $G$-exact sequence
\begin{equation}\label{00}
\begin{split}
\xymatrix{
0\ar[r]&G_{n}\ar[r] &\cdots\ar[r]&G_{0}\ar[r]^{f_{0}} &T\ar[r] &0
}
\end{split}\tag{$\sharp$}
\end{equation}
with $G_{i}$ being Gorenstein projective module for any $0\leq i\leq n$. Set $\mathrm{Ker} f_{0}=K_{1}$.
Since $i^{\ast}$ is exact, from the sequence (\ref{00}), we can obtain a new exact sequence
$$\xymatrix{
0\ar[r]&i^{\ast}G_{n}\ar[r] &\cdots\ar[r]&i^{\ast}G_{0}\ar[r] &i^{\ast}T\ar[r] &0
}$$
By the theorem \ref{result1} (2), we know that $i^{\ast}G_{i}$ is Gorenstein projective in $\mathscr A'$ for any $0\leq i\leq n$.
So the Gorenstein projective dimension of $i^{\ast}T$ is finite.
By the theorem 2.10 in \cite{Holm}, we have that $\mathrm{pd}_{G}i^{\ast}T \leq n$.

$\mathrm{Gext}^{i\geq1}(i^{\ast}T,~i^{\ast}T^{(I)})\cong \mathrm{Gext}^{i\geq1}(T,~i_{\ast}i^{\ast}T^{(I)})=0$ since $i_{\ast}i^{\ast}T\subseteq T^{G\bot}$.
Note $i^{\ast}$ is exact and preserves Gorenstein projective object.
Let $A^{\bullet}$ be a left Gorenstein projective resolution of $T$, then $i^{\ast}A^{\bullet}$ is
a left Gorenstein projective resolution of $i^{\ast}T$. By definition of the functor $\mathrm{Gext}^{i\geq1}(-,~-)$,
the isomorphism above holds. Consequently, $i^{\ast}T$ is partial $n$-Gorenstein tilting.

(2) Applying the functor $j^{\ast}$ to the sequence (\ref{00}), we can obtain a new exact sequence
\begin{equation}\label{01}
\begin{split}
\xymatrix{
0\ar[r]&j^{\ast}G_{n}\ar[r] &\cdots\ar[r]&j^{\ast}G_{0}\ar[r] &j^{\ast}T\ar[r] &0
}
\end{split}\tag{$\sharp1$}
\end{equation}
By the proposition \ref{pro1}, $j^{\ast}G_{i}$ is a Gorenstein projective module for any $0\leq i\leq n$
in $\mathscr A''$. For any Gorenstein projective module $M$ in $\mathscr A''$, we have the following commutative diagram.
$$\xymatrix@C=0.5cm{
0\ar[r]&\mathrm{Hom}_{\mathscr A''}(M,~j^{\ast}G_{n})\ar[d]^{\cong}\ar[r]&\cdots\ar[r]&\mathrm{Hom}_{\mathscr A''}(M,~j^{\ast}G_{0})\ar[d]^{\cong}\ar[r]&\mathrm{Hom}_{\mathscr A''}(M,~j^{\ast}T)\ar[d]^{\cong}\ar@{.>}[r]&0\\
0\ar[r]&\mathrm{Hom}_{\mathscr A'}(j_{!}M,~G_{n})\ar[r]&\cdots\ar[r]&\mathrm{Hom}_{\mathscr A'}(j_{!}M,~G_{0})\ar[r]&\mathrm{Hom}_{\mathscr A'}(j_{!}M,~T)\ar[r]&0
}$$
Since $j_{\ast}$ is exact, by the lemma \ref{lem1}, $j^{\ast}$ satisfies the condition (SP2). And then $j_{!}M$ is Gorenstein projective in $\mathscr A$ by the theorem \ref{result1}.
Thus the second row in above diagram is exact, i.e., the first row is exact.
From the sequence (\ref{01}), we know that $\mathrm{pd}_{G}j^{\ast}T \leq n$.

Note that $j^{\ast}$ is exact and preserves Gorenstein projective objects. Thus we have that
$\mathrm{Gext}^{i\geq1}(j^{\ast}T,~j^{\ast}T^{(I)})\cong \mathrm{Gext}^{i\geq1}(T,~j_{\ast}j^{\ast}T^{(I)})=0$ since $j_{\ast}j^{\ast}T\subseteq T^{G\bot}$.

For any Gorenstein projective module $M$ in $\mathscr A''$, then $j_{!}M$ is Gorenstein projective in $\mathscr A$ by the above discussion. By the definition, we have the following $G$-exact sequence
$$\xymatrix{
0\ar[r]&j_{!}M\ar[r] & T_{0}\ar[r] &\cdots\ar[r]& T_{n}\ar[r] &0
}$$
with $T_{i}\in \mathrm{Add}T$ for any $i$. Applying the exact $j^{\ast}$ to the above $G$-exact sequence, we give a new exact sequence as follows.
$$\xymatrix{
0\ar[r]&j^{\ast}j_{!}M\cong M\ar[r] & j^{\ast}T_{0}\ar[r] &\cdots\ar[r]& j^{\ast}T_{n}\ar[r] &0
}$$
We claim that the above sequence is $G$-exact. Indeed, let $N$ be any Gorenstein projective module in $\mathscr A''$, we can obtain the following commutative diagram.
$$\xymatrix@C=0.5cm{
0\ar[r]&\mathrm{Hom}_{\mathscr A''}(N,~j^{\ast}j_{!}M)\ar[d]^{\cong}\ar[r]&\mathrm{Hom}_{\mathscr A''}(N,~j^{\ast}T_{0})\ar[d]^{\cong}\ar[r]&\cdots\ar[r]&\mathrm{Hom}_{\mathscr A''}(N,~j^{\ast}T_{n})\ar[d]^{\cong}\ar@{.>}[r]&0\\
0\ar[r]&\mathrm{Hom}_{\mathscr A'}(j_{!}N,~j_{!}M)\ar[r]&\mathrm{Hom}_{\mathscr A'}(j_{!}N,~T_{0})\ar[r]&\cdots\ar[r]&\mathrm{Hom}_{\mathscr A'}(j_{!}N,~T_{n})\ar[r]&0
}$$
By the theorem \ref{result1}, $j_{!}N$ is Gorenstein projective in $\mathscr A$, thus the second row in above diagram is also exact. i.e.,
Our claim is correct.

To sum up, $j^{\ast}T$ is $n$-Gorenstein tilting in $\mathscr A''$.

\mskip\
\noindent\textbf{3.2 Gorenstein syzygy modules} Let $R$ be an ring. Denoted all right $R$-modules by $\mathrm{Mod}R$.

\bg{Def}\rm(\cite[Definition~2.3]{HH})
For a positive integer $n$, a module $A \in\mathrm{Mod}R$ is called a Gorenstein n-syzygy module
of M if there exists an exact sequence $0\longrightarrow A\longrightarrow X_{n-1}\longrightarrow \cdots \longrightarrow X_{0}\longrightarrow M\longrightarrow 0$
in $\mathrm{Mod}R$ with $X_{i}$ being  Gorenstein projective for any $0\leq i\leq n-1$.
\ed{Def}

By the theorem \ref{result1}, we can obtain the following statement.

\bg{Cor}\label{}%
In a recollement $R(\mathscr A', ~\mathscr A, ~\mathscr A'')$.
Let $A\in\mathscr A$ be Gorenstein n-syzygy of M.

$(1)$ If $j_{\ast}$ is exact, then $j^{\ast}A$ is Gorenstein n-syzygy of $j^{\ast}M$ in $\mathscr A''$;

$(2)$ If $i^{\ast}$ is weak perfect, $i^{\ast}$ is exact and $i_{\ast}$ satisfies the condition $\mathrm{(SP2)}$, then $i^{\ast}A$ is Gorenstein n-syzygy of $i^{\ast}M$ in $\mathscr A'$;
\ed{Cor}

Recall that a module $B \in\mathrm{Mod}R$ is called an $n$-syzygy module
of $N$ if there exists an exact sequence $0\longrightarrow B\longrightarrow P_{n-1}\longrightarrow \cdots \longrightarrow P_{0}\longrightarrow N\longrightarrow 0$
with $P_{i}$ being projective for any $0\leq i\leq n-1$.
By the proposition \ref{pro1}, we easily prove that the following result hold.

\bg{Cor}\label{g1}%
In a recollement $R(\mathscr A', ~\mathscr A, ~\mathscr A'')$.
Let $A\in\mathscr A$ be an n-syzygy module of M.

$(1)$ If $j_{\ast}$ is exact, then $j^{\ast}A$ is an n-syzygy module of $j^{\ast}M$ in $\mathscr A''$;

$(2)$ If $i^{\ast}$ is exact, then $i^{\ast}A$ is an n-syzygy module of $i^{\ast}M$ in $\mathscr A'$;
\ed{Cor}

Note that a module $A$ in $\mathrm{Mod}R$ is an $n$-syzygy module (of $M$) if and only if it is a Gorenstein $n$-syzygy module (of $N$) by \cite[Theorem 2.4]{HH}.
In general, $M$ and $N$ are different. Thus we have the following conclusion by Corollary \ref{g1}.

\bg{Cor}\label{g2}%
In a recollement $R(\mathscr A', ~\mathscr A, ~\mathscr A'')$.
Let $A\in\mathscr A$ be a Gorenstein n-syzygy module of M.

$(1)$ If $j_{\ast}$ is exact, then $j^{\ast}A$ is a Gorenstein n-syzygy module of $N$ in $\mathscr A''$;

$(2)$ If $i^{\ast}$ is exact, then $i^{\ast}A$ is a Gorenstein n-syzygy module of $L$ in $\mathscr A'$;
\ed{Cor}

Note that $N$ (resp. $L$) and $j^{\ast}M$ (resp. $i^{\ast}M$) in Corollary \ref{g2} are different in general.

\section{Resolving subcategories}

In this subsection, firstly, we introduced the definition of the relative projective dimensions and global dimensions.
Secondly, for a recollement $R(\mathscr A', ~\mathscr A, ~\mathscr A'')$ of Abelian categories, we show that a resolving subcategory in $\mathscr A$ induce
two resolving subcategories in $\mathscr A'$ and $\mathscr A''$. On the other hand, we also prove that two resolving subcategories in $\mathscr A'$ and
$\mathscr A''$ can induce a resolving subcategory in $\mathscr A$.
Finally, we give the size relationship between the relative global dimensions of $\mathscr A'$, $\mathscr A$ and $\mathscr A''$.

\bg{Def}\label{}
Let $\mathscr X$ be a resolving subcategory of $\mathscr A$ and $M\in \mathscr A$. The $\mathscr X$-projective dimension of $M$,
denoted by $pd_{\mathscr X}M$, is defined as

$\mathrm{pd}_{\mathscr X}M$ :=$inf\{n\geq0\mid$ there exists an exact sequence
$$\xymatrix{
0\ar[r] &P_{n}\ar[r] &P_{n-1}\ar[r]&\cdots\ar[r]&P_{0}\ar[r]&M\ar[r]&0
}$$
in $\mathscr A$ with $P_{i}\in \mathrm{Proj}(\mathscr A)$ for any $0\leq i\leq n\}$.

The $\mathscr X$-global dimension of $\mathscr A$, denoted by $\mathrm{gl.dim}_{\mathscr X} \mathscr A$, is defined as
$$\mathrm{gl.dim}_{\mathscr X} \mathscr A~:= \mathrm{sup}\{\mathrm{pd}_{\mathscr X}M|M\in \mathscr A\}$$
\ed{Def}

The following two lemmas are very useful. The proof of the second lemma is easy, which is omitted here.

\bg{Lem}\rm(\cite[Lemma~3.4]{ZX})\label{lem3}
Let $\mathscr X$ be a resolving subcategory of $\mathscr A$ and $0\longrightarrow L\longrightarrow M\longrightarrow N\longrightarrow0$
be an exact sequence in $\mathscr A$.

$(1)$ If $N\in \mathscr X$, then $\mathrm{pd}_{\mathscr X}L=\mathrm{pd}_{\mathscr X}M$.

$(2)$ If $M\in \mathscr X$, then $\mathrm{pd}_{\mathscr X}N=\mathrm{pd}_{\mathscr X}L+1$.

$(3)$ If $L\in \mathscr X$, then $\mathrm{pd}_{\mathscr X}M=\mathrm{pd}_{\mathscr X}N$.
\ed{Lem}

\bg{Lem}\label{lem4}
Let $\mathscr X$ be a resolving subcategory of $\mathscr A$ and $0\longrightarrow A\longrightarrow B\longrightarrow C\longrightarrow0$
be an exact sequence in $\mathscr A$. Then the following conclusions hold.

$(1)$ $\mathrm{pd}_{\mathscr X}B\leq \mathrm{max}\{\mathrm{pd}_{\mathscr X}A, \mathrm{pd}_{\mathscr X}C\}$.

$(2)$ $\mathrm{pd}_{\mathscr X}C\leq \mathrm{max}\{\mathrm{pd}_{\mathscr X}A+1, \mathrm{pd}_{\mathscr X}B\}$.

$(3)$ $\mathrm{pd}_{\mathscr X}A\leq \mathrm{max}\{\mathrm{pd}_{\mathscr X}B, \mathrm{pd}_{\mathscr X}C-1\}$.
\ed{Lem}

Next, we will the first results of this section, see Proposition \ref{pro2} and \ref{pro3}.

\bg{Pro}\label{pro2}%
Let $R(\mathscr A', ~\mathscr A, ~\mathscr A'')$ be a recollement of Abelian categories and
$\mathscr X$ be a resolving subcategory of $\mathscr A$.

$(1)$ If $i_{\ast}i^{\ast}\mathscr X\subseteq \mathscr X$, then $i^{\ast}\mathscr X$ is resolving in $\mathscr A'$.

$(2)$ If $j_{\ast}$ is exact and $j_{!}j^{\ast}\mathscr X\subseteq \mathscr X$ or $j_{\ast}j^{\ast}\mathscr X\subseteq \mathscr X$,
then $j^{\ast}\mathscr X$ is resolving in $\mathscr A''$.
\ed{Pro}

\Pf. (1) (a) For any $X\in \mathscr X$, if $i^{\ast}X=M\bigoplus N$, then $i_{\ast}M\bigoplus i_{\ast}N=i_{\ast}i^{\ast}X\in \mathscr X$.
Since $\mathscr X$ is closed under direct summands, $i_{\ast}M\in \mathscr X$, and then we have that $M\cong i^{\ast}i_{\ast}M\in i^{\ast}\mathscr X$
by the proposition \ref{propR}. i.e., $i^{\ast}\mathscr X$ is closed under direct summands.

(b) By the lemma \ref{lem1} and (a), it is clear that $\mathrm{Proj}\mathscr A'\subseteq i^{\ast}\mathscr X$.

(c) For any exact sequence $\xymatrix{0\ar[r] &i^{\ast}L\ar[r] &M\ar[r] &i^{\ast}N\ar[r]&0}$ with $L$ and $N \in \mathscr X$ in $\mathscr A'$.
Then we can get a new exact sequence  $\xymatrix{0\ar[r] &i_{\ast}i^{\ast}L\ar[r] &i_{\ast}M\ar[r] &i_{\ast}i^{\ast}N\ar[r]&0}$ in $\mathscr A$
since $i_{\ast}$ is exact. Note that both $i_{\ast}i^{\ast}L\cong L$ and $i_{\ast}i^{\ast}N\cong N$ are in $\mathscr X$ by the proposition \ref{propR}.
It follows from $\mathscr X$ is closed under extensions that $i_{\ast}M\in\mathscr X$. Thus $M\cong i^{\ast}i_{\ast}M\in i^{\ast}\mathscr X$.
i.e., $i^{\ast}\mathscr X$ is closed under extensions.

(d) For any exact sequence $\xymatrix{0\ar[r] &L\ar[r] &i^{\ast}M\ar[r] &i^{\ast}N\ar[r]&0}$ with $M$ and $N \in \mathscr X$ in $\mathscr A'$.
It is easy to verify that $L \in \mathscr X$, similar to (c). Consequently, $i^{\ast}\mathscr X$ is resolving in $\mathscr A'$.

(2) The proof is similar to (1).
\ \hfill $\Box$

\bg{Pro}\label{pro3}%
Let $R(\mathscr A', ~\mathscr A, ~\mathscr A'')$ be a recollement of Abelian categories.
Assume that $\mathscr X'$ and $\mathscr X''$ are resolving in $\mathscr A'$ and $\mathscr A''$, respectively.
If both $i^{\ast}$ and $j_{\ast}$ are exact, then the following subcategory

\centerline{$\mathscr X:=\{X\in \mathscr A|~i^{\ast}X\in \mathscr X'$ and $j^{\ast}X\in \mathscr X''\}$}
\noindent is a resolving subcategory of $\mathscr A$.
\ed{Pro}

\Pf. (a) For any $X=M\bigoplus N\in \mathscr X$, then $i^{\ast}M\bigoplus i^{\ast}N=i^{\ast}X\in \mathscr X'$.
So we have that $i^{\ast}M\in \mathscr X'$ since $\mathscr X'$ is resolving. Similarly, $j^{\ast}M\bigoplus j^{\ast}N=j^{\ast}X\in \mathscr X''$.
and then $j^{\ast}M\in \mathscr X''$ since $\mathscr X''$ is resolving. Thus $M\in\mathscr X$ by the definition.
i.e., $\mathscr X$ is closed under direct summands.

(b) By the lemma \ref{lem1}, it is easy to verify that $\mathrm{Proj}(\mathscr A)\subseteq \mathscr X$.

(c) For any exact sequence $\xymatrix{0\ar[r] &L\ar[r] &M\ar[r] &N\ar[r]&0}$ with $L$ and $N \in \mathscr X$ in $\mathscr A$.
Then we can get a new exact sequence $\xymatrix{0\ar[r] &i^{\ast}L\ar[r] &i^{\ast}M\ar[r] &i^{\ast}N\ar[r]&0}$ in $\mathscr A'$
since $i^{\ast}$ is exact. Note that both $i^{\ast}L$ and $i^{\ast}N$ are in $\mathscr X'$ by the definition.
It follows from $\mathscr X'$ is closed under extensions that $i^{\ast}M\in\mathscr X'$. Similarly, we can prove that $j^{\ast}M\in\mathscr X''$
Thus $M\in\mathscr X$. i.e., $\mathscr X$ is closed under extensions.

(d) For any exact sequence $\xymatrix{0\ar[r] &L\ar[r] &M\ar[r] &N\ar[r]&0}$ with $M$ and $N \in \mathscr X$ in $\mathscr A$.
It is easy to verify that $L \in \mathscr X$, similar to (c). Consequently, $\mathscr X$ is resolving in $\mathscr A$.
\ \hfill $\Box$

\mskip\

From the proposition \ref{pro2} and \ref{pro3}, we can obtain the following conclusion directly.

\bg{Cor}\label{}%
Let $R(\mathscr A', ~\mathscr A, ~\mathscr A'')$ be a recollement of Abelian categories.

$(1)$ If $i_{\ast}i^{\ast}\mathcal{GP}(\mathscr A)\subseteq \mathcal{GP}(\mathscr A)$, then $i^{\ast}\mathcal{GP}(\mathscr A)$ is resolving.

$(2)$ If $j_{\ast}$ is exact and $j_{!}j^{\ast}\mathcal{GP}(\mathscr A)\subseteq \mathcal{GP}(\mathscr A)$ or $j_{\ast}j^{\ast}\mathcal{GP}(\mathscr A)\subseteq \mathcal{GP}(\mathscr A)$,
then $j^{\ast}\mathcal{GP}(\mathscr A)$ is resolving.

$(3)$ If both $i^{\ast}$ and $j_{\ast}$ are exact, then the following subcategory

\centerline{$\mathscr X:=\{X\in \mathscr A|~i^{\ast}X\in \mathcal{GP}(\mathscr A')$ and $j^{\ast}X\in \mathcal{GP}(\mathscr A'')\}$}
\noindent is a resolving subcategory of $\mathscr A$.
\ed{Cor}

In order to give the second conclusion of this section, we need the following notion.
Let $R(\mathscr A', ~\mathscr A, ~\mathscr A'')$ be a recollement of Abelian categories.
and $\mathscr X$ be a resolving subcategory of $\mathscr A$. Denote the $\mathscr X$-global dimension of $\mathscr A'$ as follows.

\centerline{gl. dim$_{\mathscr X}\mathscr A'$=sup$\{\mathrm{pd}_{\mathscr X}i_{\ast}A'|A'\in \mathscr A'\}.$}
\noindent Obviously, gl. dim$_{\mathscr X}\mathscr A'\leq$ gl. dim$_{\mathscr X}\mathscr A$.

\bg{Lem}\label{lem5}
Let $R(\mathscr A', ~\mathscr A, ~\mathscr A'')$ be a recollement of Abelian categories and
$\mathscr X$ be a resolving subcategory of $\mathscr A$. If $j_{\ast}$ is exact and $j_{!}j^{\ast}\mathscr X\subseteq \mathscr X$,
then for any $A''\in\mathscr A''$, we have that the following inequations hold.

$(1)$ $\mathrm{pd}_{\mathscr X}j_{!}A''\leq \mathrm{pd}_{j^{\ast}\mathscr X}A''+\mathrm{gl. dim}_{\mathscr X}\mathscr A'+1.$

$(2)$ $\mathrm{gl. dim}_{\mathscr X}\mathscr A\leq \mathrm{gl. dim}_{\mathscr X}\mathscr A'+\mathrm{gl. dim}_{j^{\ast}\mathscr X}\mathscr A''+1.$
\ed{Lem}

\Pf. Note that $j^{\ast}\mathscr X$ is resolving in $\mathscr A''$ by Proposition \ref{pro2}.

(1) Set $\mathrm{pd}_{j^{\ast}\mathscr X}A''=s$ and
$\mathrm{gl. dim}$$_{\mathscr X}\mathscr A'=t$ are finite. We will use induction on $s$ to prove this inequation.
If $s=0$, then $A''\in j^{\ast}\mathscr X$. By the assumption, $j_{!}A''\in j_{!}j^{\ast}\mathscr X\subseteq\mathscr X$, and then $\mathrm{pd}_{\mathscr X}j_{!}A''=0.$
Thus the above inequation holds.

Assumption that the above inequation holds for less than $s-1$ ($s\geq1$). Consider the following exact sequence
$$\xymatrix{
0\ar[r] &j^{\ast}X_{s}\ar[r] &j^{\ast}X_{s-1}\ar[r]&\cdots\ar[r]&j^{\ast}X_{0}\ar[r]&A''\ar[r] &0
}$$
with $X_{i}\in \mathscr X$ for any $0\leq i\leq s$, thus, we can obtain the following exact sequence
$$(\ddagger1):\xymatrix{0\ar[r] &A''_{1}\ar[r] &j^{\ast}X_{0}\ar[r] &A''\ar[r] &0}$$
with $\mathrm{pd}_{j^{\ast}\mathscr X}A''_{1}=s-1$. Applying right exact functor $j_{!}$ to $(\ddagger1)$, we have the following two exact sequences
$$(\ddagger2):\xymatrix{0\ar[r] &Y\ar[r] &j_{!}j^{\ast}X_{0}\ar[r] &j_{!}A''\ar[r] &0}$$
and
$$(\ddagger3):\xymatrix{0\ar[r] &Y_{1}\ar[r] &j_{!}A''_{1}\ar[r] &Y\ar[r] &0.}$$
We claim that $j^{\ast}Y_{1}=0$. In fact, applying exact functor $j^{\ast}$ to the exact sequence
$$\xymatrix{0\ar[r]&Y_{1}\ar[r] &j_{!}A''_{1}\ar[r] &j_{!}j^{\ast}X_{0}\ar[r] &j_{!}A''\ar[r] &0,}$$
we have the following commutative diagram:
$$\xymatrix{
0\ar[r]&j^{\ast}Y_{1}\ar[r]\ar[d]&j^{\ast}j_{!}A''_{1}\ar[d]^{\cong}\ar[r] &j^{\ast}j_{!}j^{\ast}X_{0}\ar[d]^{\cong}\ar[r] &j^{\ast}j_{!}A''\ar[d]^{\cong}\ar[r] &0\\
&0\ar[r]&A''_{1}\ar[r] &X_{0}\ar[r] &A''\ar[r] &0
}$$
So $j^{\ast}Y_{1}=0$, and then there is an object $Y_{2}\in \mathscr A'$ such that $i_{\ast}Y_{2}=Y_{1}$.
Then $\mathrm{pd}_{\mathscr X}Y_{1}=\mathrm{pd}_{\mathscr X}i_{\ast}Y_{2}\leq \mathrm{gl. dim}$ $_{\mathscr X}\mathscr A'= t$.
Using inductive assumption for $A''_{1}$, we have the following inequation holds.
\mskip\
\centerline{$\mathrm{pd}_{\mathscr X}j_{!}A''_{1}\leq \mathrm{pd}_{j^{\ast}\mathscr X}A''_{1}$+$\mathrm{gl. dim}$$_{\mathscr X}\mathscr X$+1=$s-1+t+1=s+t$.}
\mskip\

\noindent For the sequence $(\ddagger3)$, by Lemma \ref{lem4}, we have that  $\mathrm{pd}_{\mathscr X}Y\leq s+t$.
For the sequence $(\ddagger2)$, by Lemma \ref{lem3}, we have that  $\mathrm{pd}_{\mathscr X}j_{!}A''=\mathrm{pd}_{\mathscr X}Y+1\leq s+t+1$.

(2) Set $\mathrm{gl. dim}_{\mathscr X}\mathscr A'=m<\infty$, $\mathrm{gl. dim}_{j^{\ast}\mathscr X}\mathscr A''=n<\infty$.
For any $A\in\mathscr A$, by the proposition \ref{propR}, there is an exact sequence
$$\xymatrix{0\ar[r]&i_{\ast}A'\ar[r]&j_{!}j^{\ast}A\ar[r]&A\ar[r]^{f}&i_{\ast}i^{\ast}A\ar[r]&0\\
}$$
with $A'\in\mathscr A'$ and $K=\mathrm{Ker}f$. By the lemma \ref{lem4} and (1), we have that
\begin{align*}
\mathrm{pd}_{\mathscr X} A & \leq  \mathrm{max}\{\mathrm{pd}_{\mathscr X}K, \mathrm{pd}_{\mathscr X}i_{\ast}i^{\ast}A\}\\
& \leq  \mathrm{max}\{\mathrm{pd}_{\mathscr X}i_{\ast}A'+1, \mathrm{pd}_{\mathscr X}j_{!}j^{\ast}A, m\}\\
& \leq  \mathrm{max}\{\mathrm{pd}_{\mathscr X}i_{\ast}A'+1, \mathrm{pd}_{j^{\ast}\mathscr X}j^{\ast}A+\mathrm{gl. dim}_{\mathscr X}\mathscr A'+1, m\}\\
& \leq  \mathrm{max}\{m+1, n+m+1, m\}=m+n+1
\end{align*}
So we completed the proof.
\ \hfill $\Box$

\bg{Lem}\label{lem6}
Let $R(\mathscr A', ~\mathscr A, ~\mathscr A'')$ be a recollement of Abelian categories and
$\mathscr X$ be a resolving subcategory of $\mathscr A$.
If $i_{\ast}i^{\ast}\mathscr X\subseteq \mathscr X$, then we have that $\mathrm{gl. dim}_{\mathscr X}\mathscr A'\leq\mathrm{gl. dim}_{i^{\ast}\mathscr X}\mathscr A'$.
\ed{Lem}

\Pf. By the proposition \ref{pro2}, we have that $i^{\ast}\mathscr X$ is resolving.
Set $\mathrm{gl. dim}_{i^{\ast}\mathscr X}\mathscr A'=a$. For any $A'\in\mathscr A'$, if $A'\in i^{\ast}\mathscr X$, then
$\mathrm{pd}_{\mathscr X}i_{\ast}A'=0\leq a$ by the assumption. If $A'$ is not in $i^{\ast}\mathscr X$, set $\mathrm{pd}_{i^{\ast}\mathscr X}A'=b\leq a$,
then there is an exact sequence
$$\xymatrix{
0\ar[r] &Y_{b}\ar[r] &Y_{b-1}\ar[r]&\cdots\ar[r]&Y_{0}\ar[r]&A'\ar[r]&0
}$$
with $Y_{i}\in i^{\ast}\mathscr X$ for any $0\leq i\leq b$. Applying the exact functor $i_{\ast}$ to the above exact sequence,
we can obtain the following exact sequence
$$\xymatrix{
0\ar[r] &i_{\ast}Y_{b}\ar[r] &i_{\ast}Y_{b-1}\ar[r]&\cdots\ar[r]&i_{\ast}Y_{0}\ar[r]&i_{\ast}A'\ar[r]&0
}$$
Since $i_{\ast}i^{\ast}\mathscr X\subseteq \mathscr X$, $i_{\ast}Y_{i}\in\mathscr X$ for any $i$.
So we have that $\mathrm{pd}_{\mathscr X}i_{\ast}A'\leq b\leq a$. Consequently, we can get that
$\mathrm{gl. dim}_{\mathscr X}\mathscr A'\leq\mathrm{gl. dim}_{i^{\ast}\mathscr X}\mathscr A'$.
\ \hfill $\Box$

\mskip\

Now, we can give part of the theorem \ref{th13}, which is an generalization of \cite[Theorem~4.1]{PC} in the relative homological angebra.
For another part, look at Proposition \ref{prop4}.

\bg{Th}\label{th2}%
Let $R(\mathscr A', ~\mathscr A, ~\mathscr A'')$ be a recollement of Abelian categories and
$\mathscr X$ be a resolving subcategory of $\mathscr A$. If $i_{\ast}i^{\ast}\mathscr X\subseteq \mathscr X$,
$j_{!}j^{\ast}\mathscr X\subseteq \mathscr X$ and $j_{\ast}$ is exact, then the following inequations hold.

$(1)$ $\mathrm{gl. dim}_{\mathscr X}\mathscr A\leq \mathrm{gl. dim}_{i^{\ast}\mathscr X}\mathscr A'+\mathrm{gl. dim}_{j^{\ast}\mathscr X}\mathscr A''+1.$

$(2)$ $\mathrm{gl. dim}_{j^{\ast}\mathscr X}\mathscr A''\leq \mathrm{gl. dim}_{\mathscr X}\mathscr A$.
\ed{Th}

\Pf. (1) By the lemma \ref{lem5} and lemma \ref{lem6}, the result is clear.

(2) Suppose that $\mathrm{gl. dim}_{\mathscr X}\mathscr A=p<\infty$. For any $A''\in\mathscr A''$, $j_{!}A''\in\mathscr A$. Set $\mathrm{pd}_{\mathcal{X}}j_{!}A''=q\leq p$.
Thus there is an exact sequence
$$\xymatrix{
0\ar[r]&Y_{q}\ar[r] &Y_{q-1}\ar[r]&\cdots\ar[r]&Y_{0}\ar[r]&j_{!}A''\ar[r]&0
}$$
with $Y_{i}\in \mathscr X$ for any $0\leq i\leq q$. Applying the exact functor $j^{\ast}$ to the above exact sequence,
we can obtain the following exact sequence
$$\xymatrix{
0\ar[r]&j^{\ast}Y_{q}\ar[r] &j^{\ast}Y_{q-1}\ar[r]&\cdots\ar[r]&j^{\ast}Y_{0}\ar[r]&j^{\ast}j_{!}A''\cong A''\ar[r]&0
}$$
Note that $j^{\ast}Y_{i}\in j^{\ast}\mathcal{X}$ for any $0\leq i\leq q$. i.e., $\mathrm{pd}_{j^{\ast}\mathcal{X}}A''\leq q\leq p.$
Thus we have that $\mathrm{gl. dim}_{j^{\ast}\mathscr X}\mathscr A''\leq \mathrm{gl. dim}_{\mathscr X}\mathscr A$.
\ \hfill $\Box$

\mskip\
From the theorem \ref{th2}, we can know the size relationship between the global dimensions of $\mathcal{A'}$, $\mathcal{A}$ and $\mathcal{A''}$.
We naturally want to ask who is bigger between the global dimensions of $\mathcal{A'}$ and $\mathcal{A}$. Consequently, we have the following results.

\bg{Pro}\label{prop4}
Let $R(\mathscr A', ~\mathscr A, ~\mathscr A'')$ be a recollement of Abelian categories and
$\mathscr X$ be a resolving subcategory of $\mathscr A$. If $i_{\ast}i^{\ast}\mathscr X\subseteq \mathscr X$ and $i^{\ast}$ is exact,
then $\mathrm{gl. dim}_{i^{\ast}\mathscr X}\mathscr A'\leq\mathrm{gl. dim}_{\mathscr X}\mathscr A$.
\ed{Pro}

\Pf. We know that $i^{\ast}\mathscr X$ is resolving by the proposition \ref{pro2}.
Set $\mathrm{gl. dim}_{\mathscr X}\mathscr A=m<\infty$. For any $A'\in\mathscr A'$, suppose that $\mathrm{pd}_{i^{\ast}\mathscr X}A'=n$.
Then there is an exact sequence
$$\xymatrix{
0\ar[r]&Y_{n}\ar[r] &Y_{n-1}\ar[r]&\cdots\ar[r]&Y_{0}\ar[r]&A'\ar[r]&0
}$$
with $Y_{i}\in i^{\ast}\mathscr X$ for any $0\leq i\leq n$. Applying the exact functor $i_{\ast}$ to the above exact sequence,
we can obtain the following exact sequence
$$\xymatrix{
0\ar[r]&i_{\ast}Y_{n}\ar[r] &i_{\ast}Y_{n-1}\ar[r]&\cdots\ar[r]&i_{\ast}Y_{0}\ar[r]&i_{\ast}A'\ar[r]&0
}$$
Since $i_{\ast}i^{\ast}\mathscr X\subseteq \mathscr X$, $i_{\ast}Y_{n}\in\mathscr X$ for any $0\leq i\leq n$. From the above sequence, we have that
$\mathrm{pd}_{\mathscr X}i_{\ast}A'\leq n$. We claim that $\mathrm{pd}_{\mathscr X}i_{\ast}A'=n$.
In fact, if $\mathrm{pd}_{\mathscr X}i_{\ast}A'=n-1$, Then there is an exact sequence
$$\xymatrix{
0\ar[r]&X_{n-1}\ar[r] &X_{n-2}\ar[r]&\cdots\ar[r]&X_{0}\ar[r]&i_{\ast}A'\ar[r]&0
}$$
with $X_{i}\in \mathscr X$ for any $0\leq i\leq n-1$. Applying the exact functor $i^{\ast}$ to the above exact sequence,
we can obtain the following exact sequence
$$\xymatrix{
0\ar[r]&i^{\ast}X_{n-1}\ar[r] &i^{\ast}X_{n-2}\ar[r]&\cdots\ar[r]&i^{\ast}X_{0}\ar[r]&i^{\ast}i_{\ast}A'\cong A'\ar[r]&0
}$$
So we gent that $\mathrm{pd}_{i^{\ast}\mathscr X}A'\leq n-1$. This is a contradiction, and then $\mathrm{pd}_{\mathscr X}i_{\ast}A'=n$.
Since $n=\mathrm{pd}_{\mathscr X}i_{\ast}A'\leq \mathrm{gl. dim}_{\mathscr X}\mathscr A=m$, the statement holds.
\ \hfill $\Box$

\mskip\

From the proposition \ref{prop4} and the theorem \ref{th2}, we can obtain the following two conclusions directly.

\bg{Cor}
Let $R(\mathscr A', ~\mathscr A, ~\mathscr A'')$ be a recollement of Abelian categories. If $i_{\ast}i^{\ast}\mathcal{GP}(\mathscr A)\subseteq \mathcal{GP}(\mathscr A)$,
$j_{!}j^{\ast}\mathcal{GP}(\mathscr A)\subseteq \mathcal{GP}(\mathscr A)$ and $j_{\ast}$ is exact, then the following inequations hold.

$(1)$ $\mathrm{gl. dim}_{\mathcal{GP}(\mathscr A)}\mathscr A\leq \mathrm{gl. dim}_{i^{\ast}\mathcal{GP}(\mathscr A)}\mathscr A'+\mathrm{gl. dim}_{j^{\ast}\mathcal{GP}(\mathscr A)}\mathscr A''+1.$

$(2)$ $\mathrm{gl. dim}_{j^{\ast}\mathcal{GP}(\mathscr A)}\mathscr A''\leq \mathrm{gl. dim}_{\mathcal{GP}(\mathscr A)}\mathscr A$.
\ed{Cor}

\bg{Cor}
Let $R(\mathscr A', ~\mathscr A, ~\mathscr A'')$ be a recollement of Abelian categories.
If $i_{\ast}i^{\ast}\mathcal{GP(A)}\subseteq \mathcal{GP(A)}$ and $i^{\ast}$ is exact,
then $\mathrm{gl. dim}_{i^{\ast}\mathcal{GP(A)}}\mathscr A'\leq\mathrm{gl. dim}_{\mathcal{GP(A)}}\mathscr A$.
\ed{Cor}

{\small

}

\end{document}